\documentclass{amsart}

\usepackage{amsmath,amssymb}
\usepackage{mathrsfs}
\usepackage{enumerate}

\usepackage[english, francais]{babel}
\newtheorem{theorem}{Theorem}[section]

\newtheorem{e-proposition}[theorem]{Proposition}

\newtheorem{e-definition}[theorem]{Definition\rm}
\newtheorem{conjecture}{Conjecture}

\newtheorem{theoreme}{Th\'eor\`eme}[section]
\newtheorem{lemme}[theoreme]{Lemme}

\newtheorem{frenchconjecture}{Conjecture}

\usepackage[running]{lineno} 
\linenumbers 

\def\og{\leavevmode\raise.3ex\hbox{$\scriptscriptstyle\langle\!\langle$~}}
\def\fg{\leavevmode\raise.3ex\hbox{~$\!\scriptscriptstyle\,\rangle\!\rangle$}}

\title{D\'ecomposition monomorphe des structures relationnelles et profil de classes h\'er\'editaires}
\thanks{Les auteurs d\'edient ce texte \`a la m\'emoire de Roland Fra\"{i}ss\'e et Claude Frasnay.}
\author[D.Oudrar] {Djamila Oudrar}
\address{Faculty of Mathematics, USTHB, Algiers, Algeria}
\email {dabchiche@usthb.dz; djoudrar@gmail.com}
\author [M.Pouzet]{Maurice Pouzet}
\address{ICJ, Math\'ematiques, Universit\'e
Claude-Bernard Lyon1, 43 Bd. 11 Novembre 1918
F$69622$ Villeurbanne  cedex, France and University of Calgary, Department of Mathematics and Statistics, Calgary, Alberta, Canada T2N 1N4} \email{
pouzet@univ-lyon1.fr }

\begin{document}

\maketitle


\medskip
\selectlanguage{francais}
\begin{center}
{\small Re\c{c}u le *****~; accept\'e apr\`es r\'evision le +++++\\
Pr\'esent\'e par Â£Â£Â£Â£Â£}
\end{center}

\begin{abstract}

\selectlanguage{francais}
Nous pr\'esentons une approche  structurelle de r\'esultats de sauts  dans le comportement du profil de classes h\'er\'editaires de structures finies. Nous partons de la notion suivante due \`a N.Thi\'ery et au second auteur.
Une \emph{d\'ecomposition monomorphe} d'une structure relationnelle  $R$ est  une partition de son domaine  $V(R)$ en une famille de parties $(V_x)_{x\in X}$ telles que les  restrictions de $R$ \`a deux parties  finies $A$ et  $A'$ de $V(R)$ sont isomorphes pourvu  que les traces  $A\cap V_x$ et  $A'\cap V_x$ aient m\^eme cardinalit\'e pour tout  $x\in X$. Soit  $\mathscr S_\mu$ la classe des structures relationnelles de signature  $\mu$ qui n'ont pas de d\'ecomposition  monomorphe finie. Nous montrons que si une sous classe h\'er\'editaire $\mathscr D$ de $\mathscr S_\mu $ est form\'ee  de structures relationnelles ordonn\'ees elle contient un ensemble fini $\mathfrak A$ tel que tout \'el\'ement  de  $\mathscr D$ abrite un \'el\'ement de $\mathfrak A$. En outre,  pour chaque $R\in \mathfrak A$, le profil de l'\^age    $\mathcal A(R)$  de $R$ (consistant en les sous-structures finies de $R$) est au moins exponentiel. Il en r\'esulte que, si  le profil d'une classe h\'er\'editaire de structures ordonn\'ees  n'est pas born\'e par un polyn\^ome, il  est au moins exponentiel.  Un r\'esultat faisant partie d'une classification  obtenue par   Balogh, Bollob\'{a}s et  Morris en 2006  pour les graphes ordonn\'es. {\it Pour citer cet article: Djamila Oudrar, Maurice Pouzet, C. R.
Acad. Sci. Paris, Ser. I.}
\vskip 0.5\baselineskip

\selectlanguage{english} \noindent{\bf Abstract} \vskip
0.5\baselineskip \noindent {\bf Monomorphic decomposition of relational structures. Application to the profile of hereditary classes.}
We present a structural approach of some results about jumps  in the  behavior of the profile (alias generating function) of hereditary classes of finite structures. We start with the following notion due to N.Thi\'ery and the second author.
A \emph{monomorphic decomposition} of a relational structure $R$ is a  partition of its domain $V(R)$ into a family of sets $(V_x)_{x\in X}$ such that the restrictions of $R$ to two finite subsets $A$ and $A'$  of $V(R)$ are isomorphic provided that the traces $A\cap V_x$ and $A'\cap V_x$ have the same size for each $x\in X$. Let  $\mathscr S_\mu $ be the class of relational structures of signature $\mu$ which do not have a finite monomorphic decomposition. We show that if a hereditary subclass $\mathscr D$ of $\mathscr S_\mu $ is made of  ordered relational structures then it contains a finite subset $\mathfrak A$ such that every member of $\mathscr D$ embeds  some member of $\mathfrak A$. Furthermore,  for each $R\in \mathfrak A$ the  profile of  the  age $\mathcal A(R)$  of $R$ (made of finite substructures of $R$) is at least exponential.   We  deduce  that if the profile of a hereditary class of finite ordered structures  is not bounded by a polynomial then it is at least exponential. This result is a part of  classification obtained by  Balogh, Bollob\'{a}s and Morris (2006) for ordered graphs. {\it To cite this article: Djamila Oudrar, Maurice Pouzet, C. R.
Acad. Sci. Paris, Ser. I.}
\end{abstract}

\selectlanguage{english}
\section*{Abridged English version}

Let us call   \textit{profile} of a class  $\mathscr {C}$  of finite relational structures  the integer function  $\varphi _{\mathscr{C}}$  which counts for each non negative integer $n$  the number of members of   $\mathscr{C}$
on $n$ elements, isomorphic structures being  identified.  Numerous papers discuss the behavior of this function
when  $\mathscr{C}$ is
\emph{hereditary} (that is contains  every substructure of a member of  $\mathscr{C}$) and is made of graphs (directed or not), of tournaments, of ordered sets, of ordered  graphs, of ordered hypergraphs. Futhermore, thanks to a result of Cameron \cite{cameron},  it turns out that the line of study about permutations (see \cite{A-A}) originating in the Stanley-Wilf conjecture, solved by   Marcus and Tardos (2004) \cite{Marcus},
falls under the frame of the profile of hereditary classes of relational structures (see \cite{oudrar-pouzet}).
The results show that the profile cannot be arbitrary: there are jumps in its possible growth rate. Typically, its growth  is polynomial or faster than every polynomial (\cite {pouzet.tr.1978} for ages, see  \cite{pouzet} for a survey) and for several classes of structures, either at least exponential (e.g. for tournaments \cite{BBM07,Bou-Pouz}, ordered graphs and hypergraphs \cite{B-B-M06,B-B-M06/2,klazar1} and permutations \cite{kaiser-klazar}) or  at least with the growth of the partition function (e.g. for graphs \cite {B-B-S-S}). For more, see  the survey of Klazar \cite{klazar}. A structural approach of  jump results, at least from polynomial growth to faster growth seems to be possible in view of the following result.
\begin{theorem}\label{thm:minimal1}
If the profile of a hereditary class $\mathscr C$ of finite relational structures (with a finite signature $\mu$) is not bounded by a polynomial then it contains a hereditary class  $\mathfrak A$ with this property which is minimal w.r.t. inclusion.
\end{theorem}

Indeed, according to Theorem \ref{thm:minimal1}, the  jumps in the profile are given by the growth rate of the profile of these minimal classes.  Trivially, these classes are up-directed w.r.t. embeddability, hence,  according to an old and well know result of Fra\"{\i}ss\'e, each one is the \emph{age} $\mathcal A(R)$ of some relational structure $R$ (the collection of  finite structures which are embeddable into $R$). A description of these   $R$ would allow   to evaluate the profile of $\mathcal A(R)$. We partially do  it when they are ordered.

Theorem \ref{thm:minimal1} appears in a somewhat equivalent form as Theorem 0.1 of \cite{P-T-2013}. It is not trivial, the main argument relies on a result going back to the thesis of the second author \cite{pouzet.tr.1978}, namely Lemma 4.1 p. 23 of  \cite{P-T-2013}.
We give below   an outline of  proof in the case of ordered structures.

Our first result is this:
 \begin{theorem}\label{thm:poly-expo1} (a) If an ordered  relational structure $R$ with a finite signature $\mu$, has a finite monomorphic decomposition then  the profile of $\mathcal {A}(R)$ is eventually a quasi polynomial, otherwise,  it is exponential.  More generally, (b) if $\mathscr C$ is a hereditary class of finite ordered relational structures with a finite signature $\mu$, then either there is  finite bound on the number of monomorphic components of each member of $\mathscr C$ and the profile is eventually a quasi-polynomial, or the profile of $\mathscr C$ is at least  exponential.
\end{theorem}

The proof of Theorem \ref{thm:minimal1} follows from (a). Indeed, let
$\mathscr C$ with a non polynomially bounded profile.  If   $\mathscr C$ has no infinite antichain then the collection of its hereditary subclasses  is well-founded. And we are done.  If it contains an infinite antichain then it contains an  age $A(R)$  with infinitely many bounds  which  contains no infinite antichain (this is standard,  see the mention made p.23, line 9 of \cite{P-T-2013}). Thanks to a result of Higman on words \cite{higman.1952},  the age of a structure  having a finite monomorphic decomposition is hereditary well quasi ordered, hence by the main result of \cite{pouzet 72} it has finitely many bounds.
Thus $R$ cannot have a  finite decomposition and hence by (a) its profile grows exponentially.
Since this age has no infinite antichain,  we are back to the first case. The proof of (b) follows from (a) by the same considerations.

 Let $\mathscr S_{\mu}$ be the class of all relational structures of signature $\mu$, $\mu$ finite,  without any finite monomorphic decomposition.
\begin{conjecture}
There is a finite subset $\mathfrak A$ made of incomparable structures of $\mathscr S_{\mu}$ such that every member of $\mathscr S_{\mu}$ embeds some member of $\mathfrak A$.
\end{conjecture}
If we replace $\mathscr S_{\mu}$ by the hereditary class $\mathscr D$ made of bichains then  $\mathfrak A$ has twenty elements \cite{mont-pou}.   If $\mathscr D$ is made of tournaments then  $\mathfrak A$ has twelve elements \cite{Bou-Pouz}. Our second result is this:
 \begin{theorem}
 (a) Conjecture 1 holds if we replace $\mathscr S_{\mu}$ by the hereditary class $\mathscr D$ made of ordered relational structures of arity $\mu$.
(b) If $\mathscr D$ is made of undirected graphs, $\mathfrak A$ has ten elements.  (c) If $\mathscr D$ is made of ordered reflexive (or irreflexive) directed graphs, $\mathfrak A$ contains  two hundred and thirteen elements such that the linear  order is isomorphic to $\omega$. In this last item $\mathfrak A$ is made of members whose profile  grows exponentially, in fact as fast as the Fibonacci sequence.
\end{theorem}

Our tools are Ramsey's theorem and
the notion of monomorphic decomposition of a relational structure. It was introduced in \cite{P-T-2013} in the sequel of  R.~Fra\"{\i}ss\'e who invented the notion of monomorphy  and C.~Frasnay who proved the central result about this notion \cite{frasnay 65}.
\selectlanguage{francais}

\selectlanguage{francais}

\section{Notions de base  et r\'esultats principaux}

Notre terminologie est celle de Fra\"{\i}ss\'{e} \cite{fraisse}. Une \textit{structure relationnelle}  est une paire $R:= (V,(\rho_i)_{i\in I})$ form\'ee d'une famille de relations $n_i$-aires $\rho_i$ sur $V$; l'ensemble $V$ est le \emph{domaine}, not\'e $V(R)$, la famille $\mu:=(n_i)_{i\in I}$ est la \emph{signature}. Une \emph{structure relationnelle binaire}, simplement \emph{structure binaire}, est form\'ee uniquement de relations binaires. Elle est \emph{ordonn\'ee} si une des ses relations $\rho_i$, disons par exemple $\rho_1$, est un ordre lin\'eaire.
Une structure relationnelle $R$ s'\textit{abrite} dans une structure relationnelle $R'$, fait not\'e $R\leq R'$, si $R$ est isomorphe \`a une sous-structure induite de $R'$. Une classe $\mathscr{C}$ de structures est dite \textit{h\'er\'editaire} si elle contient toute structure relationnelle qui s'abrite dans un membre de $\mathscr{C}$. La classe des structures relationnelles finies de signature $\mu$ est d\'esign\'ee par $\Omega_\mu$. Elle est pr\'eordonn\'ee par la relation d'abritement. Si $\mathscr{B}$ est un sous-ensemble de $\Omega_\mu$ alors $Forb(\mathscr{B})$ est la sous-classe des membres de $\Omega_\mu$ qui n'abritent aucun membre de $\mathscr{B}$. Clairement, $Forb(\mathscr{B})$ est une classe h\'er\'editaire. Pour la r\'eciproque, notons qu'une  \textit{borne} d'une classe h\'er\'editaire $\mathscr{C}$ de structures relationnelles finies est toute $R$ qui n'appartient pas \`a $\mathscr{C}$ telle que toute $R'$ qui s'abrite strictement dans $R$ appartient \`a $\mathscr{C}$. Clairement, chaque  borne de  $\mathscr{C}$ est finie  et si $\mathscr{B(C)}$ d\'esigne la collection des bornes de $\mathscr{C}$ alors $\mathscr{C}=Forb(\mathscr{B(C)})$.  Si $\mathscr{C}$ est contenue dans $\Omega_\mu$ et $\mathcal{A}$ est un ensemble ordonn\'e, nous posons $\mathscr{C}.\mathcal{A}:=\{(R,f): R\in \mathscr{C},f: V(R)\rightarrow \mathcal{A}\}$ et $(R,f)\leq (R',f')$ si il existe un abritement $h$ de $R$ dans $R'$ tel que $f(x)\leq f'(h(x))$ pour tout $x\in V(R)$. Nous rappelons que $\mathcal{A}$ est \textit{belordonn\'e (wqo)} si $\mathcal{A}$ ne contient ni anticha\^{i}ne infinie ni cha\^{i}ne infinie strictement d\'ecroissante. Nous disons que $\mathscr{C}$ est \textit{h\'er\'editairement belordonn\'ee} si $\mathscr{C}.\mathcal{A}$ est \textit{belordonn\'ee} pour tout ensemble belordonn\'e  $\mathcal{A}$. Nous rappelons que si la signature est finie, une sous-classe de $\Omega_\mu$ qui est h\'er\'editaire et h\'er\'editairement belordonn\'ee poss\`ede un nombre fini de bornes (\cite{pouzet 72}).

Soit $R$ une structure relationnelle.  Un sous-ensemble $V'$ de $V(R)$  est un \emph{bloc  monomorphe} de $R$ si pour tout entier $k$ et pour toute paire $A,~A'$ de sous-ensembles \`a $k$ \'el\'ements de $V(R)$, les structures induites par $A$ et $A'$ sont isomorphes d\`es que $A\setminus V'=A'\setminus V'$. Une \emph{d\'ecomposition monomorphe} de $R$ est une partition $\mathscr P$ de $V(R)$ en blocs monomorphes. Un bloc monomorphe qui est  maximal pour l'inclusion est une \emph{composante monomorphe} de $R$. Les composantes monomorphes de $R$ forment une d\'ecomposition monomorphe de $R$ et toute autre d\'ecomposition monomorphe de $R$ est plus fine qu'elle (Proposition 2.12 of  \cite{P-T-2013}).  Si une structure relationnelle infinie a une d\'ecomposition monomorphe finie en $k+1$ blocs, alors trivialement son profil est born\'e par un polyn\^ome de degr\'e  $k$. En fait, et c'est le r\'esultat principal de \cite{P-T-2013}, c'est \'eventuellement un quasi-polyn\^ome dont le degr\'e est le nombre de composantes monomorphes infinies moins $1$ (\`a partir d'un certain rang, c'est  une somme $a_{k}(n)n^{k}+\cdots+ a_0(n)$ dont les coefficients $a_{k}(n), \dots, a_0(n)$ sont des fonctions p\'eriodiques). La r\'eciproque  est fausse, sauf si la structure est ordonn\'ee. C'est une cons\'equence de notre premier r\'esultat:

\begin{theoreme}\label{thm:poly-expo} (a) Si une structure ordonn\'ee $R$ d'arit\'e finie $\mu$ a une d\'ecomposition monomorphe finie alors le profil de son age $A(R)$ est \'eventuellement un quasi-polyn\^ome,   autrement  ce profil est au moins exponentiel. Plus g\'en\'eralement, (b) si  $\mathscr C$ est une classe h\'er\'editaire de structures finies (et l'arit\'e $\mu$ est finie) alors  soit il existe une borne sur le nombre de composantes monomorphes de chaque membre de $\mathscr C$ et le profil est  \'eventuellement un quasi-polyn\^ome, soit le profil est au moins exponentiel.
\end{theoreme}
Le cas d'une classe h\'er\'editaire se ram\`ene au cas d'un \^age. C'est clair si  $\mathscr C$ est une union finie d'\^age. Sinon, $\mathscr C$ contient une anticha\^{\i}ne infinie et par suite un \^age belordonn\'e   ayant une infinit\'e de bornes. Cet age  ne peut \^etre  celui   d'une structure $R$ ayant une d\'ecomposition monomorphe finie. En effet, gr\^ace \`a un r\'esultat d'Higman sur les mots \cite{higman.1952},  l'\^age d'une structure ayant une d\'ecomposition monomorphe finie est h\'er\'editairement belordonn\'e,  et donc d'apr\`es le r\'esultat de \cite{pouzet 72}  mentionn\'e ci-dessus nous avons:
\begin{lemme}\label{lemme}
L'\^age d'une structure ayant une d\'ecomposition monomorphe finie a un nombre fini de bornes.
\end{lemme}
Ainsi,  d'apr\`es (a) la croissance du profil de $A(R)$ et donc celle du profil  de $\mathscr C$ est au moins exponentiel. Notons pour le (a) que $R:= (V, \leq, (\rho_i)_{i<m})$ a une d\'ecomposition monomorphe finie si et seulement si chaque $R_i:=(V, \leq, \rho_i)$ en a une, ou encore s'il existe une partition de $V$ en un nombre fini d'intervalles de $\leq$ tels que les  injections partielles qui les pr\'eservent sont des isomorphismes locaux de $R$.

Soit $\mathscr S_{\mu}$ la classe de toutes les structures relationnelles de signature $\mu$ qui n'ont pas de d\'ecomposition monomorphe finie.
\begin{frenchconjecture}
Il existe un sous-ensemble fini $\mathfrak A$ form\'e de structures incomparables de $\mathscr S_{\mu}$ tel que tout \'el\'ement de $\mathscr S_{\mu}$ abrite un \'el\'ement de $\mathfrak A$.
\end{frenchconjecture}

Si nous rempla\c{c}ons  $\mathscr S_{\mu}$ par la classe  $\mathscr D$ form\'ee des bicha\^{\i}nes, $\mathfrak A$ poss\`ede vingt \'el\'ements \cite{mont-pou},
 tandis que si $\mathscr D$ est form\'ee de tournois, $\mathfrak A$ poss\`ede douze \'el\'ements \cite{Bou-Pouz}. Notre second resultat est le suivant:
 \begin{theoreme}\label{thm:profil}
 (a) Cette conjecture est vraie  si $\mathscr S_{\mu}$ est remplac\'ee par la classe $\mathscr D$ form\'ee des structures   ordonn\'ees d'arit\'e $\mu$.
(b) Si $\mathscr D$ est form\'ee de graphes (non dirig\'es), $\mathfrak A$ poss\`ede dix \'el\'ements.
(c) Si $\mathscr D$ est form\'ee des graphes (dirig\'es) ordonn\'es r\'eflexifs (ou irr\'eflexifs), $\mathfrak A$ poss\`ede deux cent treize \'el\'ements si l'ordre lin\'eaire est isomorphe \`a $\omega$.
Dans ce dernier cas  $\mathfrak A$ est form\'e de membres dont le profil a une croissance exponentielle, au moins \'egale \`a la croissance de la suite de Fibonacci.
\end{theoreme}

La conclusion de (c) est contenue dans le th\'eor\`eme 1.1  de Balogh, Bollob\' as et Morris (2006) \cite{B-B-M06/2}.
Dans le cas (b), les dix \'el\'ements $G_i, 1\leq i\leq 10$, de $\mathfrak A$ ont le m\^eme ensemble de sommets $V(G_i):=\mathbb N\times\{0,1\}$. Soient  $A:=\mathbb N\times\{0\}$ et $B:=\mathbb N\times\{1\}$. Pour $i:=1,2,3$, les sous-ensembles $A$ et $B$ sont des ind\'ependants et une paire $\{(n,0),(m,1)\}$ pour $n,m\in\mathbb N$ est une ar\^ete de $G_1$ si $n=m$, une ar\^ete de $G_2$ si $n\leq m$ et une ar\^ete de $G_3$ si $n\neq m$. Ainsi,  $G_1$ est la somme directe d'une infinit\'e de copies de $K_2$ (le graphe complet \`a deux sommets) et $G_2$ est le demi-biparti complet de Schmerl et Trotter. Pour $4\leq i\leq 7$, l'un des sous-ensembles $A$, $B$ est une clique et l'autre un ind\'ependant. Ainsi, l'ensemble $E(G_4)$ des ar\^etes de $G_4$ est $E(G_1)\cup\{\{(n,0),(m,0)\}, n\neq m\in\mathbb N\}$,  $E(G_5):=E(G_2)\cup\{\{(n,0),(m,0)\}, n\neq m\in\mathbb N\}$, $E(G_6):=E(G_2)\cup\{\{(n,1),(m,1)\}, n\neq m\in\mathbb N\}$ et $E(G_7):=E(G_3)\cup\{\{(n,0),(m,0)\}, n\neq m\in\mathbb N\}$. Les graphes $G_i, i=8,9,10$, sont tels que les sous-ensembles $A$ et $B$ sont tous les deux des cliques avec $E(G_8)\cap E(G_1)=E(G_1)$, $E(G_9)\cap E(G_2)=E(G_2)$ et $E(G_{10})\cap E(G_3)=E(G_3)$. Le graphe $G_8$ est le dual de $G_3$, le graphe $G_7$ le dual de $G_4$ et le graphe $G_{10}$ le dual de $G_1$.

Dans le cas (a), chaque  graphe $G\in \mathfrak A$ a un ensemble de sommets  \'egal \`a $\mathbb N\times\{0,1\}$ ou \`a  $\{0\}\cup\mathbb N\times\{0,1\}$. Si $\leq$ est isomorphe \`a $\omega$, on obtient deux cent treize graphes r\'eflexifs (ou irr\'eflexifs) et si $\leq$ est isomorphe \`a $\omega^*$, on obtient deux cent treize autres graphes (les m\^emes que le cas pr\'ec\'edent, seul l'ordre $\leq$ change). Parmi ces exemples se trouvent les douze bicha\^{\i}nes $\mathcal B:=(V,\leq,\leq')$ de \cite{mont-pou} pour lesquelles $\leq$ est isomorphe \`a $\omega$ ou \`a $\omega^*$ (pour les huit autres, l'ordre $\leq$ est de type $\alpha+\beta$ avec $\alpha, \beta \in \{\omega, \omega^{*}\}$).

\section{Outils}

Nous  d\'efinissons la  partition  en composantes monomorphes d'une structure $R$ comme suit.

Soient $x$ et  $y$ deux \'el\'ements de $V(R)$ et   $F$ un sous-ensemble fini de $V(R)\setminus \{x, y\}$. Nous disons que $x$ et  $y$  sont \emph{$F$-\'equivalents} et notons ce fait $x \simeq_{F, R} y$ si les restrictions de $R$ \`a   $\{x\}\cup F$  et $\{y\}\cup F$ sont isomorphes. Pour un entier non n\'egatif $k$, posons $x \simeq_{k, R} y$ si $x \simeq_{F, R} y$ pour toute partie  $F$ \`a $k$ \'el\'ements de $V(R)\setminus \{x,y\}$.  Nous posons $x\simeq_{\leq k, R}y$ si $x\simeq_{k', R}y$ pour tout $k'\leq k$ et $x\simeq_{ R}y$ si $x \simeq_{k, R} y$ pour tout $k$. Ceci d\'efinit trois relations d'\'equivalence sur $V(R)$. Nous avons:

\begin{lemme} Les classes d'\'equivalence  de $\simeq_R$  sont les composantes monomorphes de $R$. \end{lemme}

\begin{lemme}
Les relations d'equivalence $\simeq_{k, R}$ et $\simeq_{\leq k, R}$ co\"{\i}ncident d\`es que  $\vert V(R)\vert \geq 2k+1$.  Les relations d'\'equivalence $\simeq_{\leq 6, R}$ et  $\simeq_{R}$ co\"{\i}ncident pour toute  structure binaire $R$. Si $R$ est un graphe dirig\'e, resp.  un graphe ordonn\'e, nous pouvons remplacer $6$ par $3$,  resp. par $2$. Si $T$ est un tournoi,  le nombre de classes d'\'equivalence de $\simeq_{\leq 3, T}$ est fini si le nombre de classes d'\'equivalence de $\simeq_{\leq 2, T}$ est fini.
Il existe un entier $i(m)$ tel que pour une structure ordonn\'ee d'arit\'e au plus $m$ les relations d'\'equivalence $\simeq_{\leq i(m), R}$ et  $\simeq_{R}$ co\"{\i}ncident.
\end{lemme}
La premi\`ere affirmation s'ensuit d'un r\'esultat de Gottlieb et Kantor sur les matrices d'incidence. Le cas des structures binaires d\'ecoule du r\'esultat de reconstruction d\^u \`a Lopez \cite{lopez 78}. Le cas des graphes dirig\'es a \'et\'e obtenu ind\'ependamment par Boudabbous \cite{boudabbous}. Le cas des structures ordonn\'ees d\'ecoule d'un r\'esultat d\^u \`a Ille \cite{ille 92}. En utilisant le r\'esultat de \cite{pouzet 79} on peut montrer qu'il n'existe aucun seuil pour les relations ternaires.

Dans le cas d'une structure binaire ordonn\'ee $R:=(E,\leq ,(\rho _{i})_{i\in I})$ la relation d'\'equivalence $\simeq_{R}$ a les propri\'et\'es suivantes:
\begin{lemme}
Toute  classe d'\'equivalence ayant au moins trois \'el\'ements est un intervalle de $R$ (au sens de Fra\"{\i}ss\'e) donc un intervalle de  $\leq$. R\'eciproquement, un intervalle de $
R$ n'est pas n\'ecessairement, une classe de $1$-\'equivalence, mais tout intervalle qui est contenu dans une classe de $1$-\'equivalence est contenu dans une classe d'\'equivalence.

\end{lemme}

\begin{lemme}
Si deux classes d'\'equivalence sont telles que leur r\'eunion est un intervalle de $\leq$  alors elles ne font pas parties d'une m\^eme classe de $1$-\'equivalence.
\end{lemme}
Dans le cas particulier des graphes ordonn\'es  dirig\'es r\'eflexifs (ou irr\'eflexifs), le r\'esultat principal qui permet de d\'eterminer les \'el\'ements de $\mathfrak A$ est le lemme suivant:
\begin{lemme}
Si un graphe ordonn\'e dirig\'e $ G:=(V,\leq,\rho)$ poss\`ede une infinit\'e de classes d'\'equivalence alors ou bien $V$ contient un sous-ensemble infini $A$ tel que deux sommets distincts de $A$ sont $0$-\'equivalents mais non $1$-\'equivalents, ou bien
$V$ contient deux sous-ensembles infinis disjoints $A_1$ et $A_2$ tels que deux sommets distincts de $A_i$, $i=1,2$, sont $1$-\'equivalents mais non \'equivalents et  pour tout $x,y\in A_i$, $i\neq j\in\{1,2\}$, la trace de l'intervalle $I_{\leq}(x,y)$ sur $A_j$ est non vide .
  \end{lemme}
Si une structure $R$ poss\`ede une infinit\'e de classes de $k$-\'equivalence  alors il existe une famille de fonctions $f:  \mathbb N\rightarrow V(R)$, $g_i: [\mathbb N]^2\rightarrow V(R)$ pour $i<k-1$ telles que pour tout $n< n'\in \mathbb N$,  $f(n)$ et $f(n')$ ne sont pas $\{ g_i( n,n'): i<k-1\}$-\'equivalents. La restriction de $R$ \`a la r\'eunion des images de cette collection de fonctions a une infinit\'e de classes d'\'equivalence.  Le th\'eor\`eme de Ramsey permet de trouver un sous-ensemble infini $X\subseteq \mathbb N$ sur lequel les fonctions de cette collection sont "invariantes" pour $R$ (la notion d'invariance est expliqu\'ee dans 2.2 de \cite{Bou-Pouz}). Supposons $X= \mathbb N$, soit $F$ l'application de $\mathbb N\times \{0, \dots k\}$ dans $V(R)$  d\'efinie par $F(n,0):= f(n)$ et $F(n,i+1):= g_i(n, n+1)$. Il s'av\`ere que si $R$ est ordonn\'ee, la restriction de $R$ \`a l'image de $F$ a une infinit\'e de classes d'\'equivalence et en outre son profil  est au moins exponentiel.  De ceci on d\'eduit la preuve du Th\'eor\`eme \ref{thm:poly-expo} et le (a) du Th\'eor\`eme \ref{thm:profil}. La m\^eme approche avec plus de soin et les lemmes ci-dessus conduisent aux (b) et (c) de ce r\'esultat.

\end{document}